\documentclass[a4paper]{amsart}
\usepackage{amssymb}
\usepackage{verbatim}
\usepackage[dvipdfm]{graphicx}
\theoremstyle{plain}
  \newtheorem{thm}{Theorem}[section]
  \newtheorem{lem}[thm]{Lemma}
  \newtheorem{cor}[thm]{Corollary}
  \newtheorem{prop}[thm]{Proposition}
  \newtheorem{conj}[thm]{Conjecture}
  
  \newtheorem{obs}[thm]{Observation}
  \newtheorem*{obs*}{Observation}
\theoremstyle{definition}
  \newtheorem{defn}[thm]{Definition}

\theoremstyle{remark}
  \newtheorem{rem}[thm]{Remark}
  
  \newtheorem*{ack}{Acknowledgments}
\newcommand{\Z}{\mathbb{Z}}
\newcommand{\C}{\mathbb{C}}

\newcommand{\Vol}{\operatorname{Vol}}
\newcommand{\CS}{\operatorname{CS}}

\newcommand{\Tr}{\operatorname{Tr}}
\newcommand{\Hom}{\operatorname{Hom}}
\newcommand{\Ad}{\operatorname{Ad}}
\newcommand{\Ker}{\operatorname{Ker}}
\newcommand{\T}{{\Sigma}}
\newcommand{\uHitoshi}{{u}}
\newcommand{\vHitoshi}{{v}}
\newcommand{\uSergei}{\tilde{u}}
\newcommand{\vSergei}{\tilde{v}}
\renewcommand{\Re}{\operatorname{Re}}
\renewcommand{\Im}{\operatorname{Im}}
\numberwithin{equation}{section}
\allowdisplaybreaks
\begin{document}
\title[$SL(2,\C)$ Chern-Simons theory and the colored Jones polynomial]
{$SL(2,\C)$ Chern-Simons theory and the asymptotic behavior of the colored Jones polynomial}

\author{Sergei Gukov}
\address{Department of Physics and Mathematics,
     California Institute of Technology, M/C 452-48,
     Pasadena, CA 91125, USA}
\email{gukov@theory.caltech.edu}
\author{Hitoshi Murakami}
\address{
Department of Mathematics,
Tokyo Institute of Technology,
Oh-okayama, Meguro, Tokyo 152-8551, Japan
}
\email{starshea@tky3.3web.ne.jp}
\date{\today}
\begin{abstract}
It has been proposed that the asymptotic behavior of the colored
Jones polynomial is equal to the perturbative expansion of the
Chern-Simons gauge theory with complex gauge group $SL(2,\C)$
on the hyperbolic knot complement.
In this note we make the first step toward verifying this
relation beyond the semi-classical approximation.
This requires a careful understanding of some delicate issues,
such as normalization of the colored Jones polynomial and
the choice of polarization in Chern-Simons theory.
Addressing these issues allows us to go beyond the volume
conjecture and to verify some predictions for the behavior
of the subleading terms in the asymptotic expansion of
the colored Jones polynomial.
\end{abstract}
\maketitle
\section{Introduction}
The original volume conjecture \cite{Kashaev:LETMP97,Murakami/Murakami:ACTAM12001}
is a remarkable relation between the limit of the colored Jones polynomial, $J_N(K;q)$,
of a knot $K$ and the volume of the knot complement $S^3\setminus{K}$:
\begin{conj}[Volume Conjecture]\label{conj:volumeconj}
For a knot $K$,
\begin{equation}\label{eq:vconj}
 \lim_{N\to\infty}
  \frac{\log \vert
        J_N\left(
             K;\exp\left(2\pi\sqrt{-1}/N
                   \right)
           \right) \vert}
       {N} = \frac{1}{2\pi} {\rm Vol} (S^3\setminus{K}),
\end{equation}
where ${\rm Vol} (S^3\setminus{K})$ is the simplicial volume of
the knot complement. In particular, if $K$ is hyperbolic,
then ${\rm Vol} (S^3\setminus{K})$ is the hyperbolic volume.
\end{conj}
The physical interpretation of this relation was proposed
in \cite{Gukov:COMMP2005}, where it was conjectured that the asymptotic
expansion of the colored Jones polynomial $J_N(K;q)$ in the limit
$N \to \infty$, $q \to 1$ should be equal to the partition
function of the $SL(2,\C)$ Chern-Simons theory on the knot
complement $S^3\setminus{K}$:
\begin{conj}\label{conj:jones}
As $N\to\infty$ and $k\to\infty$ with $u:=2\pi \sqrt{-1}
(\frac{N}{k}-1)$ kept fixed, the colored Jones polynomial
has the asymptotic expansion
\begin{multline}\label{eq:asymptotics}
  \log J_N(K;\exp(2\pi\sqrt{-1}/k))
  \\
  \underset{N,k\to\infty}{\sim}
  \frac{k}{\sqrt{-1}} S(u)
  +
  \frac{1}{2} \delta_K(u)\log{k}
  +
  \frac{1}{2}\log\left(\frac{T_K(u)}{2\pi^2}\right)
  +
  \sum_{n=1}^{\infty}\left(\frac{2\pi}{k}\right)^nS_{n+1}(u)
\end{multline}
where the function $S(u)$ in the first term
is the classical action of the Chern-Simons theory;
$T_K(u)$ is the Ray--Singer torsion of the knot complement twisted
by the flat connection corresponding to the representation
$\rho:\pi_1(S^3\setminus{K}) \to SL(2;\C)$ determined by $u$
(see \cite{Ray/Singer:ADVAM11971,Ray/Singer:ANNMA11973} and
\S\ref{section:beyond} below);
the number $\delta_K(u) \in \Z$ is determined by the topology of
the knot complement and the representation $\rho$;
finally, the function $S_{n}(u)$ denotes the $n$-loop contribution
\cite{Witten:COMMP1989,Axelrod/Singer:1992,BarNatan:Thesis}.
\end{conj}
\par
In particular, this physical interpretation of the volume
conjecture opens an avenue for several generalizations.
First, it suggests that, for a knot $K$, there exists a 1-parameter
family of relations like \eqref{eq:vconj} --- sometimes called
the ``generalized'' or ``parameterized'' volume conjecture ---
which relate a family of limits of the colored Jones polynomial
to the volume function $\Vol(K;u)$ on the character variety
of the knot complement, see Conjecture \ref{conj:parameterization} below.
Moreover, the interpretation via $SL(2,\C)$ Chern-Simons theory predicts
the structure of the subleading terms in the asymptotic expansion
of the colored Jones polynomial, where each term in \eqref{eq:asymptotics}
has an {\it a priori} definition and can be computed independently,
before comparing to the colored Jones polynomial.
In what follows, first we shall discuss the leading term in the
expansion \eqref{eq:asymptotics} and then return to the subleading
terms in \S\ref{section:beyond}.
\par
Notice, the leading term in the expansion \eqref{eq:asymptotics}
gives precisely the generalization of the volume conjecture
which can be stated as follows:
\begin{conj}\label{conj:parameterization}
There exists an open subset $\mathcal{O}_{K}$ of $\C$ such that for any
$u\in\mathcal{O}_{K}$ the following limit exists:
\begin{equation}\label{eq:limit}
  \lim_{N\to\infty}
  \frac{\log
        J_N\left(
             K;\exp\left(
                     \left(u+2\pi\sqrt{-1}\right)/N
                   \right)
           \right)}
       {N} = \frac{2\pi}{(u+2\pi\sqrt{-1})} S(u).
\end{equation}
\end{conj}
The precise formulation of the Conjectures \ref{conj:jones} and \ref{conj:parameterization}
requires fixing normalizations of the colored Jones polynomial and the partition function
of the Chern-Simons theory.
Moreover, in the geometric quantization approach to Chern-Simons gauge theory,
the partition function is determined by quantizing the moduli space of
flat connections which involves an extra ambiguity,
the choice {\it polarization}.\footnote{For a general introduction to geometric
quantization see \cite{Woodhouse:Quantization}.}
In principle, all these choices should be fixed by {\it a priori} arguments,
such as a non-perturbative formulation of the $SL(2,\C)$ Chern-Simons theory.
However, with the lack of such arguments one might use the fact that all
these choices do not depend on the knot $K$ and, therefore, can be fixed
once and for all by considering, say, the leading term $S(u)$ for a particular knot.
Thus, the simplest choice of polarization (that is the choice of
the symplectic potential $S(u)$) consistent with the asymptotic
behavior of the colored Jones polynomial for the figure-eight knot at $\Re (u) =0$
gives $\Im S(u) = \frac{1}{2\pi} \Vol(K;u)$, where $\Vol(K;u)$ is
the volume function
\cite{Neumann/Zagier:TOPOL85,Cooper/Culler/Gillet/Long/Shalen:INVEM1994},
\begin{equation}\label{eq:Schlafli}
d\,\Vol (K;u) = -\frac{1}{2} \left(\Re (u) d \Im (v) -\Re (v) d
\Im (u) \right).
\end{equation}
Notice, this determines $\Im S(u)$ only up to terms proportional to $\Re (u)$
(which are highly constrained, though, see \S\ref{section:polarization} below).
Unfortunately, until recently there was no example of a hyperbolic knot for
which the asymptotic behavior of the colored Jones polynomial would be
studied for $\Re (u) \ne 0$. One simple possibility is to assume that
there are no additional terms proportional to $\Re (u)$, so that
$\Im S(u) = \frac{1}{2\pi} \Vol(K;u)$. This gives the generalized
volume conjecture originally proposed in \cite{Gukov:COMMP2005}.
\par
Recently, it was realized \cite{Murakami:2006} that the generalized
volume conjecture proposed in \cite{Gukov:COMMP2005} requires a modification,
precisely by the terms proportional to $\Re (u)$.
Indeed, the explicit study of the limit \eqref{eq:limit}
for the figure-eight knot \cite{Murakami/Yokota:JREIA}
and torus knots at $\Re (u) \ne 0$ leads to the version
of the Conjecture \ref{conj:parameterization} with the function $S'(u)$,
where \cite{Murakami:2006}:
\begin{equation}\label{eq:volume}
\Im S'(u) = \frac{1}{2\pi} \Vol(K;u) +
\frac{1}{2}\Re(u)+\frac{1}{4\pi}\Re(u)\Im(v).
\end{equation}
We believe it is this generalization of the volume conjecture --- called
the parameterized volume conjecture in \cite{Murakami:2006} --- which has 
a chance of being true.
We note (and explain in more detail in \S\ref{section:polarization})
that in the geometric quantization approach to the $SL(2,\C)$ Chern-Simons theory,
both $S(u)$ and $S'(u)$ can be interpreted as the semi-classical action of
the $SL(2,\C)$ Chern-Simons theory, obtained in a different polarization.
For this reason, the physical considerations in \cite{Gukov:COMMP2005}
can not tell us whether it is $S(u)$ or $S'(u)$ which should appear
in the right-hand side of \eqref{eq:asymptotics} and \eqref{eq:limit}.
However, given the evidence in \cite{Murakami:2006,Murakami/Yokota:JREIA},
we believe it is the function $S'(u)$ which is the correct 
semi-classical action of the $SL(2,\C)$ Chern-Simons theory,
and both Conjectures \ref{conj:jones} and \ref{conj:parameterization}
should be considered with $S'(u)$ given by \eqref{eq:volume}.
We present further evidence for this in \S\ref{section:beyond}.
\par
One of the main goals of the present paper is to show that
the Conjectures \ref{conj:jones} and \ref{conj:parameterization}
with the function $S'(u)$ are consistent with the proposal
that the asymptotic expansion of the colored Jones polynomial
is equal to the loop expansion
of the partition function of the $SL(2,\C)$ Chern-Simons theory.
In particular, in \S\ref{section:polarization} we show that
the function $S'(u)$ can be interpreted as 
the semi-classical action in the $SL(2,\C)$ Chern-Simons theory,
and differs from $S(u)$ by a choice of polarization.
Note, that $S(u) = S'(u)$ for $\Re (u) =0$, which for hyperbolic knots
is expected to correspond to the case of cone-manifolds.
Once we identify the correct choice of polarization,
in \S\ref{section:beyond} we present further evidence
for the Conjecture \ref{conj:jones} by studying the subleading
terms in the asymptotic expansion \eqref{eq:asymptotics}
of the colored Jones polynomial and compare with
the expected behavior of the partition function
of the $SL(2,\C)$ Chern-Simons theory.
\par
Before we proceed, let us describe the representation
$\rho\colon\pi_1 (S^3\setminus{K}) \to SL(2,\C)$ determined by $u$.
Following \cite{Murakami:2006}, we define
\begin{equation}\label{vviau}
v_{K}(u) := 4\pi\frac{d\,S'(u)}{d\,u}-2\pi\sqrt{-1},
\end{equation}
where $S'(u)$ is defined by \eqref{eq:limit}.
Then, $\rho$ can be defined as a representation from
$\pi_1 (S^3\setminus{K})$ to $SL(2,\C)$ sending
the longitude and the meridian to the elements whose
eigenvalues are $(l,m) = (-\exp (-v_K (u)/2),\exp(u/2))$.
The representation $\rho$ defines a flat $SL(2,\C)$ bundle
over $S^3\setminus{K}$ which we denote $E_{\rho}$.
This bundle will play an important role in \S\ref{section:beyond}.
\begin{ack}
The authors would like to thank J\'er\^ome Dubois, Stavros Garoufalidis,
and Toshiaki Hattori for helpful conversations.
It is also a pleasure to thank the organizers of the conference ``Around the Volume
Conjecture'' at Columbia University in March 2006, which stimulated much of this work.
This work was supported in part by the DOE under grant
number DE-FG03-92-ER40701, in part by RFBR grant 04-02-16880,
and in part by the grant for support of scientific schools NSh-8004.2006.2 (S.G.),
and in part by Grant-in-Aid for Scientific Research (B) (15340019) (H.M.).
\end{ack}

\section{Choice of polarization}\label{section:polarization}
In this section we show that the version of the generalized volume
conjecture \ref{conj:parameterization} proposed in \cite{Murakami:2006}
with the function $S'(u)$ is consistent with the interpretation
in terms of $SL(2;\C)$ Chern-Simons theory suggested in \cite{Gukov:COMMP2005}.
In particular, we show that $S'(u)$ can be interpreted as 
the semi-classical action in the $SL(2,\C)$ Chern-Simons theory,
and differs from $S(u)$ by a choice of polarization.
To explain this in detail, let us start by fixing notations.
We use $\uSergei$ and $\vSergei$ instead of $u$ and $v$
used in \cite{Gukov:COMMP2005} respectively. The relation between
$(\uSergei,\vSergei)$ and $(\uHitoshi,\vHitoshi)$ will be
explained later.
\par
Let us recall the argument in \cite{Gukov:COMMP2005}.
Fix an oriented (not necessarily hyperbolic) knot $K$ in $S^3$.
Let $M$ be its complement $S^3\setminus{\operatorname{int}N(K)}$,
where $N(K)$ is the tubular neighborhood of $K$ and $\operatorname{int}N(K)$
is its interior.
The boundary of $M$ is denoted by $\Sigma$.
Note that $\pi_1(\Sigma)\cong\Z\times\Z$ is generated by the meridian
$\mu$ and the longitude $\lambda$, where $\mu$ bounds a disk in
$N(K)$ and is oriented so that the linking number of $K$ and $\mu$ is $-1$,
and $\lambda$ is null homologous in $M$ and parallel to $K$ in $N(K)$.
(So our orientation of the meridian here is {\em different} from usual
one in knot theory.)
\par
We consider a representation $\rho$ of $\pi_1(M)$ to $SL(2;\C)$ and denote
by $m=\exp{\vSergei}$ and $l=\exp{\uSergei}$ the eigenvalues of its images of
$\mu$ and $\lambda$ respectively.
Then the pair $(m,l)$ is a zero of the $A$-polynomial
\cite{Cooper/Culler/Gillet/Long/Shalen:INVEM1994}.
The zero locus of the $A$-polynomial defines a Lagrangian
submanifold $L$ of $\mathcal{P}$ with respect to the $2$-form $\omega$,
where $\mathcal{P}=\C^{\ast}\times\C^{\ast}$
is the representation space of $\pi_1(\Sigma)$
and $\omega$ is defined as follows \cite[\S3]{Gukov:COMMP2005}.
\begin{equation*}
  \omega:=-\frac{1}{\pi}d\uSergei\wedge d\vSergei
\end{equation*}
so that $\omega\bigm|_L=0$.
(Note that here we put $\sigma=k$ in \cite[(3.7)]{Gukov:COMMP2005},
and rescaled $\omega$ by a factor of $k$ which now explicitly appears
as the coefficient of the classical action in \eqref{eq:asymptotics}.)
If we put
\begin{equation*}
  \theta
  :=
  \frac{1}{2\pi}
  \left(
    \vSergei\,d\uSergei-\uSergei\,d\vSergei
    +d\left(\uSergei\overline{\vSergei}\right)
  \right),
\end{equation*}
we have $d\theta=\omega$ \cite[(3.26)]{Gukov:COMMP2005}.
Let $S$ be the classical Chern--Simons action corresponding to $\rho$.
Then $S$ can be obtained by integrating $\theta$ over a path on $L$, that is,
we have
\begin{equation*}
  S=\int\theta.
\end{equation*}
Since the Lagrangian submanifold $L$ is quantizable the integral above is
well-difined \cite[\S3]{Gukov:COMMP2005}.
Note that $S$ depends on the choice of $\theta$ satisfying $d\theta=\omega$.
So we can define $S$ only up to a total derivative on $\mathcal{P}$
(the choice of polarization).
\par
One possible choice of $\theta$, consistent with $d\theta=\omega$,
gives
\begin{equation}\label{eq:S}
  S
  =
  \frac{\sqrt{-1}}{2\pi}\bigl(\Vol(l,m)+2\pi^2\sqrt{-1}\CS(l,m)\bigr),
\end{equation}
where $\Vol(l,m)$ and $\CS(l,m)$ are the volume and the
Chern--Simons invariant of the representation $\rho$. According to
\eqref{eq:asymptotics}, the leading term of the $\log$ of the
$N$-colored Jones polynomial evaluated at $e^{2 \pi \sqrt{-1}/k}$
should be $-\sqrt{-1} k S$ when $N\to\infty$ and $k\to\infty$.
This leads to a generalization of the volume conjecture
\cite[(5.12)]{Gukov:COMMP2005}
\begin{equation}\label{eq:GVC}
  \lim_{N\to\infty, k\to\infty}
  \frac{\log J_N \left(K;\exp\left(2\pi\sqrt{-1}/k\right)\right)}{k}
  =
  \frac{1}{2\pi}
  \left(
    \Vol(l,m)+2\pi^2\sqrt{-1}\CS(l,m)
  \right).
\end{equation}
\par
Now, as we pointed out earlier, the dependence on the choice of
polarization is related to the choice of the $1$-form $\theta$
such that $d\theta=\omega$. In particular, we can consider the
following $1$-form:
\begin{equation*}
  \theta'
  :=
  \frac{1}{2\pi}
  \left(-2\uSergei\,d\vSergei+2\pi\sqrt{-1}d\vSergei\right).
\end{equation*}
Note that $d\theta=d\theta'=\omega$. Let $S'$ be the classical
Chern--Simons action obtained from $\theta'$ ($S':=\int\theta'$).
Then
\begin{equation*}
\begin{split}
  dS-dS'
  &=
  \theta-\theta'
  \\
  &=
  \frac{1}{2\pi}
  d
  \left(
    \uSergei\vSergei+\uSergei\overline{\vSergei}-2\pi\sqrt{-1}\vSergei
  \right)
  \\
  &=
  \frac{1}{\pi}
  d
  \left(
    \uSergei\Re(\vSergei)-\pi\sqrt{-1}\vSergei
  \right).
\end{split}
\end{equation*}
Therefore from \eqref{eq:S} we have
\begin{equation}\label{eq:Vol}
\begin{split}
  \Vol(l,m)
  &=
  2\pi
  \Im(S)
  \\
  &=
  2\pi
  \Im(S')
  +
  2\Im(\uSergei)\Re(\vSergei)-2\pi\Re(\vSergei).
\end{split}
\end{equation}
\par
Now we consider the pair $(\uHitoshi,\vHitoshi)$ used in the previous sections.
As described in \cite{Murakami/Yokota:JREIA}, $\uHitoshi$ and $\vHitoshi$ are
related to $\uSergei$ and $\vSergei$ as follows:
\begin{equation*}
\begin{cases}
  \vSergei&=\dfrac{\uHitoshi}{2},
  \\[3mm]
  \uSergei&=-\dfrac{\vHitoshi}{2}.
\end{cases}
\end{equation*}
So from \eqref{eq:Vol} we have
\begin{equation*}
  \Vol(l,m)
  =
  2\pi \Im(S')
  -
  \pi\Re(\uHitoshi)
  -
  \frac{1}{2}\Re(\uHitoshi)\Im(\vHitoshi).
\end{equation*}
Comparing with \eqref{eq:volume} and \eqref{eq:asymptotics},
$-\sqrt{-1} k S'$ gives the leading term of the $\log$ of the
$N$-colored Jones polynomial if we use $\theta'$ to define the
classical action $S'$,
\begin{equation*}
  \lim_{N\to\infty, k\to\infty}
  \frac{\log |J_N \left(K;\exp\left(2\pi\sqrt{-1}/k\right)\right)|}{k}
  =
  \frac{1}{2\pi}
    \Vol(K;u)
  +
  \frac{1}{2}\Re(\uHitoshi)
  +
  \frac{1}{4\pi}\Re(\uHitoshi)\Im(\vHitoshi)
\end{equation*}
This is precisely the version of the generalized volume conjecture
proposed in \cite{Murakami:2006}.
\section{Beyond the Leading Order}\label{section:beyond}
Now let us discuss the subleading terms in the asymptotic
expansion \eqref{eq:asymptotics} of the colored Jones polynomial.
\par
The simplest knot to consider is the unknot $U$. Since
$J_N(U;q)=[N]:=
\left(q^{N/2}-q^{-N/2}\right)/\left(q^{1/2}-q^{-1/2}\right)$, we
have
\begin{equation}\label{eq:asymptotic_unknot}
\begin{split}
  &\log J_N(U;\exp(2\pi\sqrt{-1}/k))
  \\
  &~=
  \log\sin(N\pi/k)-\log\sin(\pi/k)
  \\
  &~\sim
  \log{k}
  -
  \log{\pi}
  +
  \log\sin\bigl(u/(2\sqrt{-1})\bigr)
  +\text{terms of the order $k^{-1}$ or lower}
\end{split}
\end{equation}
for large $N$ and $k$ with fixed $N/k=u/(2\pi\sqrt{-1})+1$ and
$u\ne0$.
\par
Since the colored Jones polynomial vanishes at $u=0$, in this case
one should use the {\em reduced} colored Jones polynomial
$V_N(K;q):=J_N(K;q)/J_N(U;q)$ to study the asymptotic expansion
\eqref{eq:asymptotics}. We should note, however, that it is $J_N(K;q)$
which naturally appears in Chern-Simons theory.
\subsection{The logarithmic term}
In this subsection we study the logarithmic term in the asymptotic
expansion \eqref{eq:asymptotics}.
\par
For a knot $K$ in $S^3$, let $M_K$ be the complement of the
interior of the regular neighborhood of $K$. Denote by $\T_K$ the
boundary torus of $M_K$. Let $H^i\left(M_K;E_{\rho}\right)\cong
H^i\left(S^3\setminus{K};E_{\rho}\right)$ be the $i$-th cohomology
group of $M_K$ with coefficients in the flat $SL(2,\C)$ bundle
$E_{\rho}$ which was described in the end of \S 1. It is well
known that $H^i\left(M_K;E_{\rho}\right)$ is isomorphic to the
cohomology $H^i\left(M_K;sl_2(\C)\right)$ with coefficient the Lie
algebra $sl_2(\C)$ twisted by the adjoint action of $\rho$. We
will mainly use $H^i\left(M_K;sl_2(\C)\right)$ for calculation.
\par
Define
\begin{equation}\label{eq:delta_rep}
  \delta^{\text{rep}}_K(\rho)
  :=
  3
  +
  h^1(M_K;\rho)
  -
  h^0(M_K;\rho),
\end{equation}
where
$h^0(M_K;\rho):=\dim H^0(M_K;\rho)$ and
\begin{equation*}
  h^1(M_K;\rho)
  :=
  \dim
  \left(
    \Ker[H^1(M_K;\rho)\to H^1(\T_K;\rho')]
  \right),
\end{equation*}
where the map is induced by the inclusion $\T_K\to M_K$ and
$\rho'$ is the restriction of $\rho$ to $\pi_1(\T_K)$. From
\cite[Page 42]{Hodgson/Kerckhoff:JDIFG21998} and
\cite[D{\'e}monstration de Proposition 3.7 (Page
72)]{Porti:MAMCAU1997}, we have
\begin{equation}\label{eq:h^1}
\begin{split}
  h^1(M_K;\rho)
  &=
  \dim H^1(M_K;\rho)-\frac{1}{2}\dim H^1(\T_K;\rho')
  \\
  &=
  \begin{cases}
    \dim H^1(M_K;\rho)-3
    &\quad\text{if $\rho'$ is trivial,}
    \\
    \dim H^1(M_K;\rho)-1
    &\quad\text{if $\rho'$ is non-trivial,}
  \end{cases}
\end{split}
\end{equation}
where $\rho'$ is trivial if it sends every element of
$\pi_1(\T_K)$ to $\pm I$ with $I$ the identity matrix.
\par
Among other things, the identification of the asymptotic expansion
\eqref{eq:asymptotics} with the partition function in the
Chern-Simons theory implies (see {\it e.g.}
\cite{Freed/Gompf:CMPHAY1991,Gukov:COMMP2005, Witten:COMMP1989})
\begin{equation}\label{logterm}
  \delta_K(u)=\delta^{\text{rep}}_K(\rho),
\end{equation}
where $\rho$ is determined by $u$ as described in \S1. The first
term (equal to $3$) in \eqref{eq:delta_rep} comes from the
normalization by the partition function of $S^3$,
\begin{equation*}
\begin{split}
\log Z(S^3) & =
  \log\sqrt{2/k}\sin(\pi/k)
  \\
  & \sim
  {}-\frac{3}{2}\log{k}+\cdots.
\end{split}
\end{equation*}
\par
To evaluate \eqref{eq:delta_rep}, it is useful to note that
\begin{equation}\label{isotropy}
  h^0(K;\rho) = \dim (H_{\rho}),
\end{equation}
where $H_{\rho}$ is the isotropy group of $E_{\rho}$, the subgroup
of $SL(2,\C)$ that commutes with the holonomies of flat
connections on $E_{\rho}$. In other words,
\begin{equation*}
H_{\rho}=\{g\in SL(2;\C)\mid\text{$g\rho(\gamma)=\rho(\gamma)g$
for any $\gamma\in\pi_1(M_K)$}\}.
\end{equation*}
Another useful fact is that $h^1(M_K;\rho)$ counts the
infinitesimal deformations of $E_{\rho}$ for a fixed
representation of $\pi_1(\T_K)$
\cite{Weil:ANNMA11960,Weil:ANNMA11962,Weil:ANNMA11964}.
\par
We can show that $h^1(M_K;\rho)=0$ if $K$ is a hyperbolic knot or a torus knot
for a typical representation $\rho$.
\begin{defn}\cite[Definition~3.21]{Porti:MAMCAU1997}
Let $\gamma$ be a simple closed curve on $\T_K$ for a knot $K$.
Then a representation $\rho$ of $\pi_1(M_K)$ in $SL(2;\C)$ is
called $\gamma$-regular if
\begin{enumerate}
\item the map $H^1(M_K;\rho)\to H^1(\gamma;\rho_0)$ induced by the
inclusion $\gamma\to M_K$ is injective, and \item if
$\Tr(\rho(\pi_1(\T_K)))=\pm1$, then $\rho(\gamma)\ne\pm I$.
\end{enumerate}
\end{defn}
The following representations are known to be $\gamma$-regular
\cite[Page 83]{Porti:MAMCAU1997}.
\begin{enumerate}
\item
The holonomy representation corresponding to the complete hyperbolic structure of a hyperbolic knot is $\mu$-regular, where $\mu$ is the meridian.
\item
For a hyperbolic knot, if the Dehn surgery along a simple closed curve $\gamma$
is hyperbolic then the holonomy representation induced by $\gamma$ is $\gamma$-regular.
\item
For a torus knot, any irreducible representation is both $\mu$-regular and
$\lambda$-regular, where $\lambda$ is the longitude
\cite[Example~1]{Dubois:CANMB2006}.
\end{enumerate}
We can show the following proposition.
\begin{prop}
Let $\rho$ be a $\gamma$-regular representation of $\pi_1(M_K)$
for a knot $K$ for some simple closed curve $\gamma$ on $\T_K$. Then
$\delta^{\text{rep}}_{K}(\rho)=3-h^0(K;\rho)$.
\par
Moreover if $\rho$ is non-Abelian, then $\delta^{\text{rep}}_{K}(\rho)=3$.
\end{prop}
\begin{proof}
Let $\rho'$ be the restriction of $\rho$ to $\pi_1(\T_K)$, where
$\T_K$ is the boundary torus of the regular neighborhood of $K$.
\par
Since $\rho$ is $\gamma$-regular, the composition
\begin{equation*}
  H^1(M_K;\rho)
  \to
  H^1(\T_K;\rho')
  \to
  H^1(\gamma;\rho_0)
\end{equation*}
is injective and so is the map $H^1(M_K;\rho)\to H^1(\T_K;\rho')$.
Therefore we have $h^1(M_{K};\rho)=\dim\Ker[H^1(M_K;\rho)\to
H^1(\T_K;\rho')]=0$ and the formula follows.
\par
If $\rho$ is non-Abelian, we have $h^0(K;\rho)=0$ (see for example
\cite[Lemme~0.7 (ii)]{Porti:MAMCAU1997}).
\end{proof}
As a corollary we have
\begin{cor}
If $K$ is a hyperbolic knot or a torus knot and $\rho$ is non-Abelian representation of $\pi_1(S^3\setminus{K})$ in $SL(2;\C)$, then
$\delta_K^{\text{rep}}=3$.
\end{cor}
\subsubsection{Abelian representations}
We will study the case where $\rho$ is Abelian.
\begin{lem}
For any knot $K$, there exists open sets $U_{1}\ni I$ and $U_{2}\ni -I$ of
$SL(2;\C)$ such that if $\rho$ is an Abelian representation sending the meridian into $U_1$ or $U_2$ then $\delta^{\text{rep}}_K(\rho)=2$.
\end{lem}
\begin{proof}
First note that $H_1(M_K;\Z)\cong\Z$ is generated by the meridian.
We choose an element $\mu\in\pi_1(M_K)$ that is mapped to the meridian by the
Abelianization.
If $\rho$ is Abelian, it is determined by the image $A:=\rho(\mu)\in SL(2;\C)$,
since it factors through $H_1(M_K;\Z)$.
\par
We can calculate $H^i(M_K;\rho)$
by using the infinite cyclic covering space $\widetilde{M_K}$
of $M_K$.
We have the following chain complex of $\widetilde{M_K}$ as
$\C[t,t^{-1}]$-modules by using the Fox differential calculus
as described in the proof of \cite[Theorem~6.1]{Boden/Herald/Kirk/Klassen:GEOTO2001}:
\begin{equation*}
  C_2\xrightarrow{d_2}C_1\xrightarrow{d_1}C_0,
\end{equation*}
where $C_2$ is generated by $\{r_1,r_2,\dots,r_{n-1}\}$,
$C_1$ is generated by $\{s_1,s_2,\dots,s_n\}$,
$C_0$ is generated by $\{p\}$, and $d_2$ and $d_1$ are given as follows.
We define $d_1$ by $d_1(s_i):=(t-1)p$ for $1\le i\le n$.
Let $F(t)$ be the $n\times(n-1)$ matrix with entries in $Z[t,t^{-1}]$
given by the Fox free differential calculus.
Then $d_2(r_i)=\sum_{j=1}^{n}F_{ji}(t)s_j$, where $F_{ji}(t)$ is the $(j,i)$ entry of $F(t)$.
Note that the sum of each colums of $F(t)$ is zero
(and so $d_1\circ d_2=0$), and that the determinant of any $(n-1)\times(n-1)$ matrix obtained from $F(t)$ by deleting
any row gives the Alexander polynomial $\Delta(K;t)$ of $K$.
(See for example \cite[Chapter~11]{Lickorish:1997}.)
\par
Then the twisted cohomology $H^i(M_K;sl_2(\C))$ is calculated from
the following cochain complex:
\begin{multline}\label{eq:cochain}
  \{0\}\to
  \Hom_{\C[t,t^{-1}]}\left(C_0;sl_2(\C)\right)
  \xrightarrow{d^{\ast}_1}
  \Hom_{\C[t,t^{-1}]}\left(C_1;sl_2(\C)\right)
  \\
  \xrightarrow{d^{\ast}_2}
  \Hom_{\C[t,t^{-1}]}\left(C_2;sl_2(\C)\right).
\end{multline}
Here an element $\varphi\in\Hom_{\C[t,t^{-1}]}\left(C_0;sl_2(\C)\right)$ is given by $\varphi(p)\in sl_2(\C)$ in such a way that
\begin{equation*}
  \varphi(g(t)p):=g(\Ad(A))\varphi(p),
\end{equation*}
where $g(t)\in\C[t,t^{-1}]$ is a Laurent polynomial and $\Ad(A)$ is the adjoint representation.
Since
\begin{equation*}
  d^{\ast}_1(\varphi)(s_i)
  :=
  \varphi(d_1(s_i))
  =
  \varphi((t-1)p)
  =
  (\Ad(A)-I)\varphi(p)
  =
  A\varphi(p)A^{-1}-\varphi(p)
\end{equation*}
for any $i$, we have $\dim\Ker{d^{\ast}_1}=1$ if $A\ne\pm I$. Thus
we have $h^0(M_K;\rho)=1$. (Another way to show this is to use
\eqref{isotropy}.)
\par
Similarly, if $\psi\in\Ker d^{\ast}_2$, then
\begin{equation*}
  d^{\ast}_2(\psi)(r_j)
  :=
  \psi(d_2(r_j))
  =
  \psi
  \left(
    \sum_{k=1}^{n}F_{kj}(t)s_k
  \right)
  =
  \sum_{k=1}^{n}F_{kj}\bigl(\Ad(A)\bigr)\psi(s_k)
  =
  O
\end{equation*}
for $j=1,2\dots n-1$, where $O$ is the $2\times2$ zero matrix.
Since $\sum_{k=1}^{n}F_{kj}=0$ for any $j$, we have
\begin{equation}\label{eq:Ker_d_2}
  \sum_{k=1}^{n-1}F_{kj}\bigl(\Ad(A)\bigr)\{\psi(s_k)-\psi(s_n)\}
  =
  O.
\end{equation}
By conjugation we may assume that $A$ is of the form
$\begin{pmatrix}a&b\\0&a^{-1}\end{pmatrix}$. Note that for any
$X:=\begin{pmatrix}x&y\\z&-x\end{pmatrix}\in sl_2(\C)$ the
$(2,1)$-entry of $A^m X A^{-m}$ is equal to $a^{-2m}z$. So
comparing the $(2,1)$-entries of \eqref{eq:Ker_d_2} we have
\begin{equation*}
  \sum_{k=1}^{n-1}F_{kj}\left(a^{-2}\right)\{\psi_{2,1}(s_k)-\psi_{2,1}(s_n)\}=0
\end{equation*}
for $j=1,2,\dots,n-1$, where $\psi_{2,1}(s_k)$ is the $(2,1)$
entry of $\psi(s_k)$. Now we know the determinant of the
$(n-1)\times(n-1)$ matrix $\left(F_{kj}(t)\right)_{1\le k,j\le
n-1}$ is $\Delta(K;t)$. Since $\Delta(K;\pm1)$ is odd for any knot
$K$ (see for example \cite[Chapter~6]{Lickorish:1997}), there
exist open sets $U_1\ni I$ and $U_2\ni -I$ in $SL(2;\C)$ such that
$\Delta(K;a^{-2})\ne0$ if $a$ is an eigenvalue of a matrix in
$U_1\cup U_2$. Thus if $A\in U_1\cup U_2$, the matrix
$\left(F_{kj}(t)\right)_{1\le k,j\le n-1}$ is non-singular and we
have $\psi_{2,1}(s_k)=\psi_{2,1}(s_n)$ for any $1\le k\le n-1$
{}from \eqref{eq:Ker_d_2}. Similarly, we have
$\psi_{1,1}(s_k)=\psi_{1,1}(s_n)$ and
$\psi_{1,2}(s_k)=\psi_{1,2}(s_n)$, which implies
$\psi(s_k)=\psi(s_n)$. This means that the $1$-cocycle group $Z^1$
of \eqref{eq:cochain} is given by
\begin{equation*}
  Z^1
  =
  \left\{
    \psi\in\Hom_{\C[t,t^{-1}]}(C_1;sl_2(\C))
    \mid
    \psi(s_1)=\psi(s_2)=\dots=\psi(s_n)
  \right\}
\end{equation*}
and so $\dim Z^1=3$.
Since $\dim\Ker{d_1^{\ast}}=1$, the dimension of $1$-coboundary group $B^1$
is $2$.
We finally have $h^1(M_K;\rho)=\dim Z^1-\dim B^1=1$.
\par
{}From \eqref{eq:h^1} we have $\delta^{\text{rep}}_K(\rho)=2$.
\end{proof}
In \cite{Garoufalidis/Le:aMMR}, S.~Garoufalidis and T.~Le prove that for any
knot $K$ the limit
\begin{equation*}
  \lim_{N\to\infty}V_N\bigl(K;\exp((u+2\pi\sqrt{-1})/N)\bigr)
\end{equation*}
exists if $u$ is sufficiently close to $-2\pi\sqrt{-1}$, which was
first proved for the figure-eight knot by the second author in
\cite{Murakami:2005}. This means that $\delta_K(u)=2$ from
\eqref{eq:asymptotic_unknot}. We expect that such $u$ determines
an Abelian representation $\rho$.
\subsubsection{Connected-sums}
We will discuss the behavior of $\delta_K(u)$ and
$\delta^{\text{rep}}_K(\rho)$ under connected-sum and satallite.
Let us denote by $K_1\sharp K_2$ the connected-sum of two knots
$K_1$ and $K_2$. Then we have $J_N(K_1\sharp
K_2;q)=J_N(K_1;q)J_N(K_2;q)/J_N(U;q)$. Therefore we have
\begin{align}\label{eq:delta_connected_sum}
  \delta_{K_1\sharp K_2}(u)
  &=
  \delta_{K_1}(u)+\delta_{K_2}(u)-2
\end{align}
if $u\ne0$. Correspondingly, for $u=0$ we have $V_N(K_1\sharp
K_2;q)=V_N(K_1;q)V_N(K_2;q)$ which implies
\begin{equation}\label{eq:delta_connected_sum_trivial}
  \delta_{K_1\sharp K_2}(0)
  =
  \delta_{K_1}(0)+\delta_{K_2}(0)
\end{equation}
\par
Now, let us compare this with $\delta^{\text{rep}}_{K_1\sharp
K_2}(\rho).$ Note that the complement $M_{K_1\sharp K_2}$ is
obtained from $M_{K_1}$ and $M_{K_2}$ by glueing along an annulus.
More precisely, $M_{K_1\sharp K_2}$ is obtained from $M_{K_1}$ and
$M_{K_2}$ by identifying annuli $A_1\in T_{K_1}$ and $A_2\in
T_{K_2}$, where $A_1$ (resp. $A_2$) is the regular neighborhood of the meridian $\mu_1$ (resp. $\mu_2$) of $K_1$ (resp. $K_2$).
\par
We first calculate $H^0(S^1;\rho)$.
\begin{lem}\label{lem:cohomology_annulus}
If a representation $\rho\colon\pi_1(S^1)\to SL(2;\C)$ is not $\pm I$
then we have
\begin{equation*}
  H^{0}(S^1;\rho)=\C,
\end{equation*}
Otherwise, we have
\begin{equation*}
  H^{0}(S^1;\rho)=\C^3,
\end{equation*}
\end{lem}
\begin{proof}
We use the interpretation of $h^0$ described in \eqref{isotropy}.
If $\rho$ is not $\pm I$ its isotropy group is one-dimensional and so
$h^0(S^1;\rho)=1$.
The case where $\rho=\pm I$ the equality follows since the cohomology group is the usual one.
\end{proof}
Next we calculate $H^i(M_{K_1\sharp K_2};\rho)$ for $i=0$ and $1$
under some assumption.
\par
\begin{lem}\label{lem:delta_connected_sum}
Let
$\rho\colon\pi_1(M_{K_1\sharp K_2})\to SL(2;\C)$
be a representation, and
$\rho_1$, $\rho_2$, and $\rho_0$ restrictions of $\rho$ to
$\pi_1(M_{K_1})$, $\pi_1(M_{K_2})$, and $\pi_1(A)$ respectively, where $A:=M_{K_1}\cap M_{K_2}(\cong A_1\cong A_2)$.
Suppose that $\rho_1$ and $\rho_2$ are $\mu_1$- and $\mu_2$-regular respectively.
\par
Then
\begin{equation}\label{eq:delta_rep_connected_sum}
  \delta^{\text{rep}}_{K_1\sharp K_2}(\rho)
  =
  \delta^{\text{rep}}_{K_1}(\rho_1)+\delta^{\text{rep}}_{K_2}(\rho_2)
  -2
\end{equation}
if $\rho_0\ne\pm I$, and
\begin{equation}\label{eq:delta_rep_connected_sum_trivial}
  \delta^{\text{rep}}_{K_1\sharp K_2}(\rho)
  =
  \delta^{\text{rep}}_{K_1}(\rho_1)+\delta^{\text{rep}}_{K_2}(\rho_2)
\end{equation}
if $\rho_0=\pm I$.
\end{lem}
\begin{proof}
The Mayer--Vietoris exact sequence for $M_{K_1\sharp K_2}=M_{K_1}\cup M_{K_2}$,
$M_{K_1}$, $M_{K_2}$ and $A$ gives the
following exact sequence:
\begin{equation*}
\begin{split}
  \{0\}&\to
  H^{0}(M_{K_1\sharp K_2};\rho)
  \to
  H^{0}(M_{K_1};\rho_1)\oplus H^{0}(M_{K_2};\rho_2)
  \to
  H^{0}(A;\rho_0)
  \\
  &\to
  H^{1}(M_{K_1\sharp M_2};\rho)
  \to
  H^{1}(M_{K_1};\rho_1)\oplus H^{1}(M_{K_2};\rho_2)
  \xrightarrow{j^{\ast}_1-j^{\ast}_2}
  H^{1}(A;\rho_0),
\end{split}
\end{equation*}
where $j^{\ast}_1$ and $j^{\ast}_2$ are induced by the inclusions
$j_1\colon A\to M_{K_1}$ and $j_2\colon A\to M_{K_2}$ respectively.
Therefore we have
\begin{equation*}
\begin{split}
  &\dim H^0(A;\rho_0)
  \\
  &~=
  \left\{
    \dim H^0(M_{K_1};\rho_1)
    +
    \dim H^0(M_{K_2};\rho_2)
    -
    \dim H^0(M_{K_1\sharp K_2};\rho)
  \right\}
  \\
  &~+
  \left\{
    \dim H^1(M_{K_1\sharp K_2};\rho)
    -
    \dim\Ker(j^{\ast}_1-j^{\ast}_2)
  \right\}.
\end{split}
\end{equation*}
\par
If $\rho_0\ne\pm I$, then from Lemma~\ref{lem:cohomology_annulus} and \eqref{eq:h^1} we have
\begin{equation*}
\begin{split}
  \delta^{\text{rep}}_{K_1\sharp K_2}
  &=
  3
  +\dim\Ker\left[H^1(M_{K_1\sharp K_2};\rho)\to
                H^1(T_{K_1\sharp K_2};\rho')
           \right]
  -\dim H^0(M_{K_1\sharp K_2};\rho)
  \\
  &=
  3
  +\dim H^1(M_{K_1\sharp K_2};\rho)-1
  -\dim H^0(M_{K_1\sharp K_2};\rho)
  \\
  &=
  2
  +
  \dim H^0(A;\rho_0)
  +
  \dim\Ker(j^{\ast}_1-j^{\ast}_2)
  -h^0(M_{K_1};\rho_1)
  -h^0(M_{K_2};\rho_2)
  \\
  &=
  3
  +
  \dim\Ker(j^{\ast}_1-j^{\ast}_2)
\end{split}
\end{equation*}
Now since $j^{\ast}_1$ and $j^{\ast}_2$ are injective and $\dim
H^1(A;\rho_0)=1$, we conclude that
$\dim\Ker(j^{\ast}_1-j^{\ast}_2)
=\dim(\Im{j^{\ast}_1}\cap\Im{j^{\ast}_2})=1$. Therefore
\begin{equation*}
  \delta^{\text{rep}}_{K_1\sharp K_2}
  =
  4
  =
  \delta^{\text{rep}}_{K_1}(\rho_1)+\delta^{\text{rep}}_{K_2}(\rho_2)
  -2.
\end{equation*}
\par
If $\rho=\pm I$, then since $\dim H^0(A;\rho_0)=3$ we have
\begin{equation*}
\begin{split}
  \delta^{\text{rep}}_{K_1\sharp K_2}
  &=
  \dim H^1(M_{K_1\sharp K_2};\rho)
  -
  \dim H^0(M_{K_1\sharp K_2};\rho)
  \\
  &=
  3
  +
  \dim\Ker(j^{\ast}_1-j^{\ast}_2)
  -h^0(M_{K_1};\rho_1)
  -h^0(M_{K_2};\rho_2).
\end{split}
\end{equation*}
But in this case we have
\begin{equation*}
\begin{split}
  &\dim\Ker(j^{\ast}_1-j^{\ast}_2)
  \\
  =&
  3
  +
  \dim\Ker[H^1(M_{K_1};\rho_1)\to H_1(T_{K_1})]
  +
  \dim\Ker[H^1(M_{K_2};\rho_2)\to H_1(T_{K_2})]
\end{split}
\end{equation*}
and so
\begin{equation*}
\begin{split}
  \delta^{\text{rep}}_{K_1\sharp K_2}
  &=
  6
  +
  h^1(M_{K_1},\rho_1)
  +
  h^1(M_{K_2};\rho_2)
  -
  h^0(M_{K_1},\rho_1)
  +
  h^0(M_{K_2};\rho_2)
  \\
  &=
  \delta^{\text{rep}}_{K_1}(\rho_1)+\delta^{\text{rep}}_{K_2}(\rho_2).
\end{split}
\end{equation*}
\end{proof}
\begin{obs}
Equation~\eqref{eq:delta_connected_sum} should be compared
with its ``physics'' counterpart, Equation~\eqref{eq:delta_rep_connected_sum}.
On the other hand, Equations~\eqref{eq:delta_connected_sum_trivial}
and \eqref{eq:delta_rep_connected_sum_trivial} which appear to
have the same form actually have very different origin.
\end{obs}
\begin{rem}
The assumption of Lemma~\ref{lem:delta_connected_sum} is true for hyperbolic knots and for torus knots.
In particular, it is true for hyperbolic knots with representations close to
the holonomy representations corresponds to their complete
hyperbolic structures.
\end{rem}
\subsubsection{Satellite knots}
We consider satellite knots.
\par
Let $K$ be a knot in a solid torus $D$.
If $e\colon D\to S^3$ is an embedding, then the image $e(K)$ forms a knot.
We call $e(K)$ a satellite of the knot $C$ with companion $K$,
where $C$ is the image of the core of $D$.
Then the complement $M_{e(K)}$ of the interior of the regular neighborhood of
$K$ is obtained from
$S^3\setminus{\operatorname{Int}{D}}$ and
$D_{K}:=D\setminus\operatorname{Int}N(K)$ by pasting along their
boundaries $\partial M_{C}$
and $\partial{D}$, where $\operatorname{Int}$ denotes the interior, $N(K)$
is the regular neighborhood of $K$ in $D$, and
$M_C:=S^3\setminus\operatorname{Int}N(C)$.
Note that $\partial D_{K}$ consists of two tori; $\partial D$ and
$\partial M_{e(K)}$.
\par
We can compute $\delta^{\text{rep}}_{e(K)}(\rho)$ under some
assumptions.
\begin{lem}\label{lem:delta_satellite}
Let $e(K)$ be a satellite of a knot $C$ with companion $K$ and
$\rho\colon\pi_1(S^3\setminus{e(K)})\to SL(2;\C)$ a non-trivial representation.
We assume the following four conditions:
\begin{gather}
  h^0(M_{C};\rho_1)=0,
  \tag{i}
  \\
  h^1(M_{C};\rho_1)=0,
  \tag{ii}
  \\
  h^0(D_{K};\rho_2)=0,
  \tag{iii}
  \\
  \Ker\left[H^1(D_{K};\rho_2)\to H^1(\partial D_{K};{\rho'}_2)\right]=0,
  \tag{iv}
\end{gather}
where $\rho_1$, $\rho_2$, and ${\rho'}_2$ are the restrictions of
$\rho$ to
$\pi_1(S^3\setminus{C})$, $\pi_1(D\setminus{K})$ and
$\partial D_{K}$ respectively, and the map in {\rm(iv)} is induced by the includion.
We also assume all the induced representations $\rho_0$, $\rho_1$, and $\rho_2$
are non-trivial.
\par
Then we have $\delta^{\text{rep}}_K(\rho)=4$.
\end{lem}
\begin{proof}
{}From the Mayer--Vietoris exact sequence for
$M_{e(K)}=M_{C}\cup D_{K}$,
$M_{C}$, $D_{K}$ and
$T=M_{C}\cap D_{K}$, we have
\begin{equation*}
\begin{split}
  \{0\}
  &\to
  H^0(M_{e(K)};\rho)
  \to
  H^0(M_{C};\rho_1)
  \oplus
  H^0(D_{K};\rho_2)
  \to
  H^0(T;\rho_0)
  \\
  &\to
  H^1(M_{e(K)};\rho)
  \to
  H^1(M_{C};\rho_1)
  \oplus
  H^1(D_{K};\rho_2)
  \xrightarrow{j^{\ast}_1-j^{\ast}_2}
  H^1(T;\rho_0),
\end{split}
\end{equation*}
where $\rho_0$ is the restriction of $\rho$ to $\pi_1(T)$,
and $j_1$ and $j_2$ are inclusions.
{}From the assumptions {\rm (i)} and {\rm (iii)}, we have
$h^0(M_{e(K)};\rho)=\dim H^0(M_{e(K)};\rho)=0$.
\par
We note that $h^1(M_{C};\rho)=\dim\Ker{j^{\ast}_1}=0$.
We also note that
\begin{equation*}
\begin{split}
  \dim\Ker{j^{\ast}_2}
  &=
  \dim\Ker
  \left[
    H^1(D_{K};\rho_2)\to
    H^1(\partial D_{K};{\rho'}_2)\to
    H^1(T;\rho_0)
  \right]
  \\
  &=
  \dim\Ker
  \left[
    H^1(\partial D_{K};{\rho'}_2)\to
    H^1(T;\rho_0)
  \right]
  \\
  &=
  1
\end{split}
\end{equation*}
since $\dim H^1(T;\rho_0)=1$ from \cite[D{\'e}monstration de
Proposition 3.7 (Page 72)]{Porti:MAMCAU1997}. Therefore the kernel
of the map $j^{\ast}_1-j^{\ast}_2$ is one-dimensional, and we have
$H^1(M_{e(K)};\rho)=H^0(T;\rho_0)\oplus\C$. Since $\pi_1(T)$ is
abelian and $\rho_0\ne\pm I$, we have $H^0(T;\rho_0)=\C$ from
\eqref{isotropy}. From \eqref{eq:h^1}, we finally get
$h^1(M_{e(K)};\rho)=1$.
\par
Therefore $\delta_{e(K)}(u)=4$ as stated.
\end{proof}
\begin{obs}
This can be compared with a result of Zheng \cite[Theorem 1.4]{Zheng:Whitehead},
where he studies Whitehead doubles of non-trivial torus knots and proves that
$\delta_K(0)=4$ for such knots.
\par
In this case $D_{K}$ is homeomorphic to the complement of
(the interior of the regular neighborhood of) the Whitehead link, which
is hyperbolic.
So from \cite[Corollary 1.2]{Hodgson/Kerckhoff:JDIFG21998},
the assumption {\rm(iv)} is satisfied if $\rho$ is a small deformation of the
holonomy representation.
The other conditions also hold.
\end{obs}
\begin{rem}
The assumption
$h^1(S^3\setminus{C};\rho_1)=h^0(S^3\setminus{C};\rho_1)=0$
of Lemma~\ref{lem:delta_satellite}
is true if $C$ is a hyperbolic knot and $\rho_1$ is close to the
holonomy representation corresponding to the complete hyperbolic structure.
The assumption
$h^1(D\setminus{K};\rho_2)
=h^0(D\setminus{K};\rho_2)=0$
may be true if $D\setminus{K}$ possesses a complete hyperbolic structure
and $\rho_2$ is close to the holonomy representation.
\end{rem}
\begin{rem}
It would be interesting to study the asymptotic behavior of the colored Jones
polynomial of the Hopf link $H$.
Since $J_N(H;q)=[N^2]$ (see for example \cite[Lemma~14.2]{Lickorish:1997}),
we have
\begin{equation*}
\begin{split}
  \log J_N(H;\exp{2\pi\sqrt{-1}/k})
  &=
  \log\sin\left(\pi\left(\frac{N}{k}\right)^2 k\right)
  -
  \log\sin(\pi/k)
  \\
  &\sim
  2\log{k}
  +
  2\log(N/k)
\end{split}
\end{equation*}
if $N/k$ is {\em very small}.
This formula would mean that $\delta_H(u)=4$, which is consistant with the fact that $3+h^1(T\times I)-h^0(T\times I)=4$ since the complement of the Hopf link is $T\times I$.
\end{rem}
\subsection{Ray--Singer torsion}
The next interesting term in \eqref{eq:asymptotics} is the term
containing the Ray--Singer torsion $T_K(u) := T(S^3\setminus{K},\rho)$.
Following the conventions used in the literature on Chern-Simons
theory, we define the Ray--Singer torsion of a 3-manifold $M$
with respect to a flat bundle $E_{\rho}$ corresponding to $\rho$ by
\begin{equation}\label{rstorsion}
  T(M,\rho)
  =
  \exp
  \left(
    - \frac{1}{2} \sum_{n=0}^3 n (-1)^n \log {\det}' \Delta^{E_\rho}_n
  \right),
\end{equation}
where $\Delta^{E_\rho}_n$ is the Laplacians on $n$-forms with coefficients in $E_{\rho}$,
and ${\det}' \Delta^{E_\rho}_n$ is the regularized determinant of the restriction
of the orthocomplement of its kernel. Using Poincar\'e duality, one finds
\begin{equation}
  T (M,\rho)
  =
  \frac{\left({\det}' \Delta_0^{E_\rho}\right)^{3/2}}
       {\left({\det}' \Delta_1^{E_\rho}\right)^{1/2}}.
\end{equation}
We remind that the relation between $u$ and the corresponding
representation $\rho$ was discussed in the end of \S 1.
Note, that the definition of the Ray--Singer torsion is particularly simple
when
$h^0(S^3\setminus{K};\rho)=h^1(S^3\setminus{K};\rho)=0$;
in this case the Laplacians $\Delta^{E_\rho}_n$ have empty kernels.
\par
Very much like the leading term in the expansion \eqref{eq:asymptotics},
$T_K(u)$ is a non-trivial function on the character variety.
In view of the Cheeger-M\"uller theorem \cite{Cheeger:ANMAAH1979,Muller:ADMTA41978},
it would be interesting to compare the Ray--Singer torsion
as a function on the character variety
to the Reidemeister torsion studied by Porti
\cite{Porti:MAMCAU1997}.
For example, for the figure-eight knot $E$ and $\Re (u)=0$ one has
\begin{align}
  T_E(u)
  &=
  \frac{1}{\sqrt{(3/2 - \cos \alpha)(1/2 + \cos \alpha) }}
  \label{eq:torsion_nontrivial}
  \\
  \intertext{if $u\ne0$, and}
  T_E(0)
  &=
  \frac{\pi^2}{\sqrt{(3/2 - \cos \alpha)(1/2 + \cos \alpha) }}
  \label{eq:torsion_trivial}
\end{align}
where $\alpha = |\sqrt{-1} u|$ is the singular (cone)
angle of the cone manifold $M = S^3 \setminus K$.
Note that in \cite[\S5.3, Exemple~1]{Porti:MAMCAU1997}
the torsion is given as $\pm1/T_E(u)$ for $0<\sqrt{-1}u<2\pi/3$.
Note also that we use the {\em reduced} colored Jones polynomial $V_N$
when $u=0$.
The difference between \eqref{eq:torsion_nontrivial} and
\eqref{eq:torsion_trivial} comes from $-\log\pi$ in
\eqref{eq:asymptotic_unknot}.
\par
Let us consider the following function of $N$ and $r:=N/k=1+u/(2\pi\sqrt{-1})$.
\begin{multline*}
  \Re
  \left\{
     \log J_N\bigl(E;\exp(2\pi{r}\sqrt{-1}/N)\bigr)
  \right\}
  \\
  -
  \Re
  \left\{
    \frac{N}{r\sqrt{-1}}S(r)
    +\frac{3}{2}\log\left(\frac{N}{r}\right)
    +\frac{1}{2}\log
     \left(
        \frac{T_E(u)}{2\pi^2}
      \right)
  \right\},
\end{multline*}
for the figure-eight knot $E$.
Here we use the result of \cite{Murakami:KYUMJ2004}.
We expect that it vanishes when $N\to\infty$.
In Figure~1--6 we use MATHEMATICA to plot the graphs of this function
for $N=100$, $200$, $300$, $400$, $500$ and $1000$.
\begin{figure}[h]
\begin{center}
\includegraphics{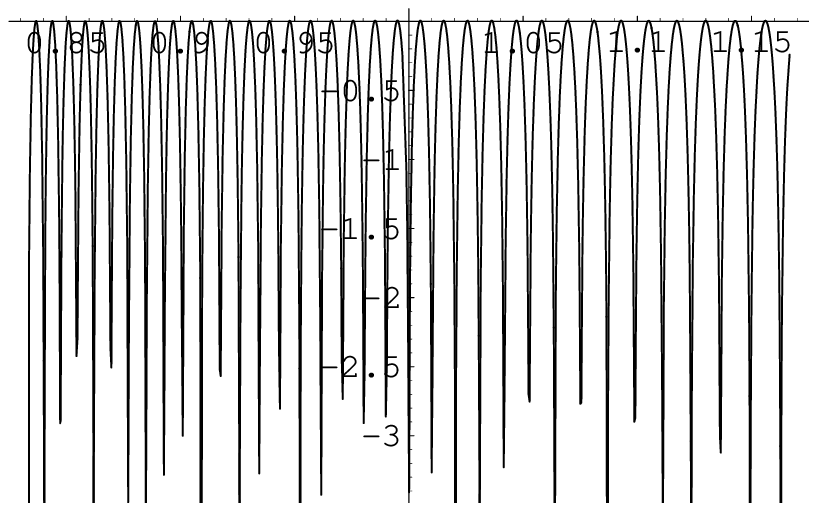}
\end{center}
\caption{$N=100$}
\end{figure}

\begin{figure}[h]
\begin{center}
\includegraphics{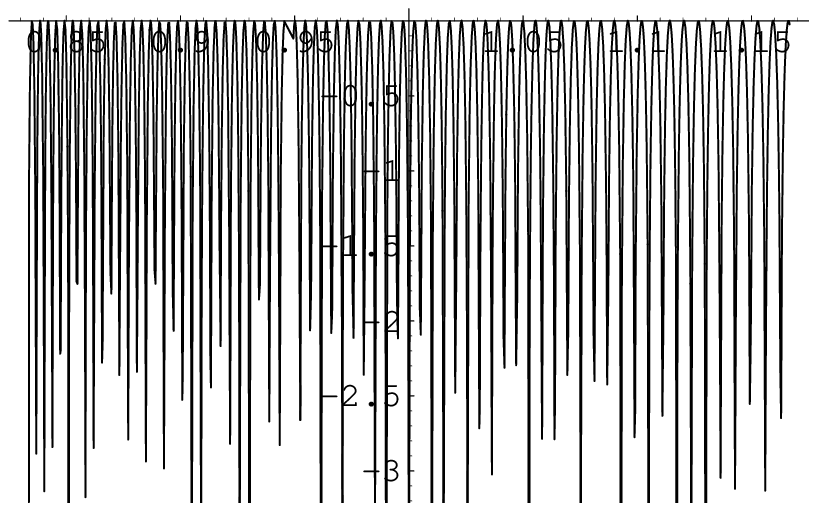}
\end{center}
\caption{$N=200$}
\end{figure}

\begin{figure}[h]
\begin{center}
\includegraphics{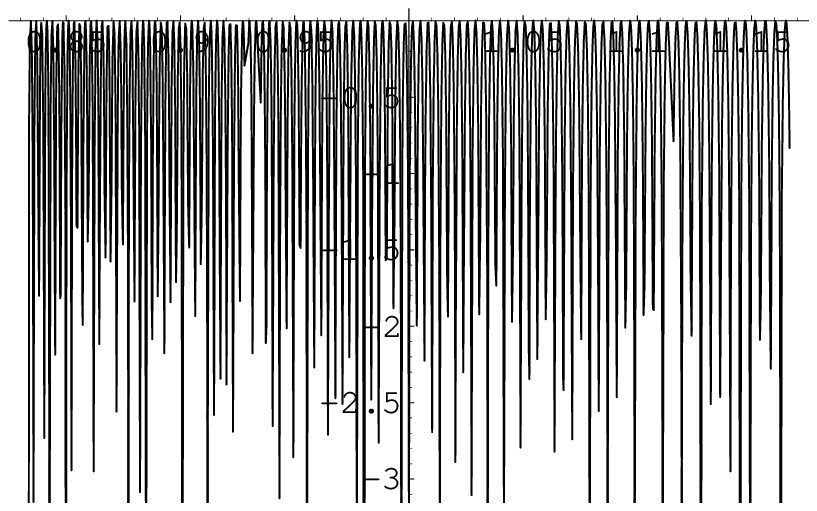}
\end{center}
\caption{$N=300$}
\end{figure}

\begin{figure}[h]
\begin{center}
\includegraphics{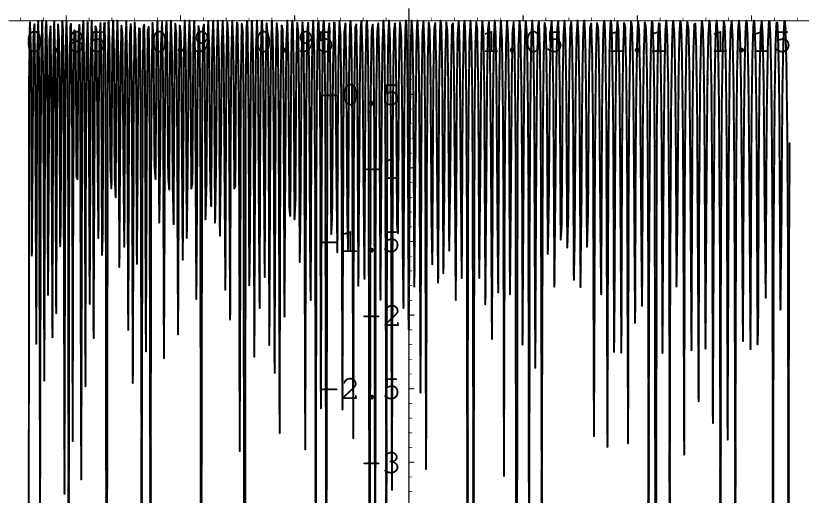}
\end{center}
\caption{$N=400$}
\end{figure}

\begin{figure}[h]
\begin{center}
\includegraphics{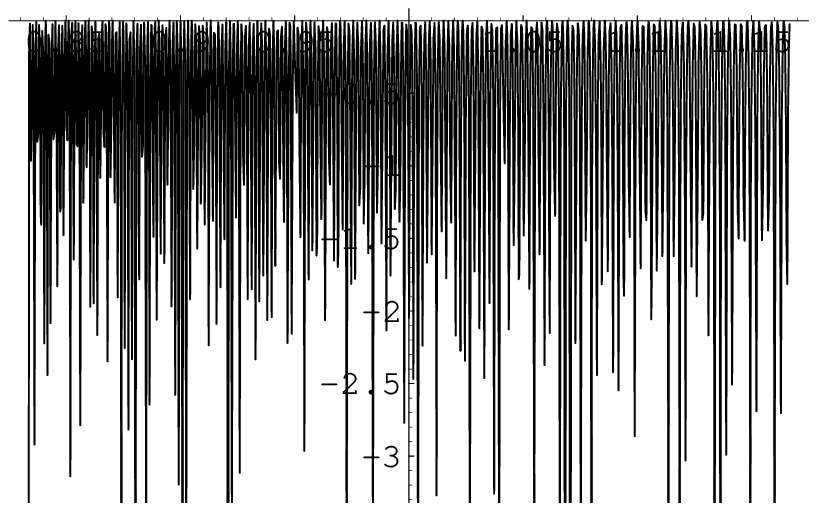}
\end{center}
\caption{$N=500$}
\end{figure}

\begin{figure}[h]
\begin{center}
\includegraphics{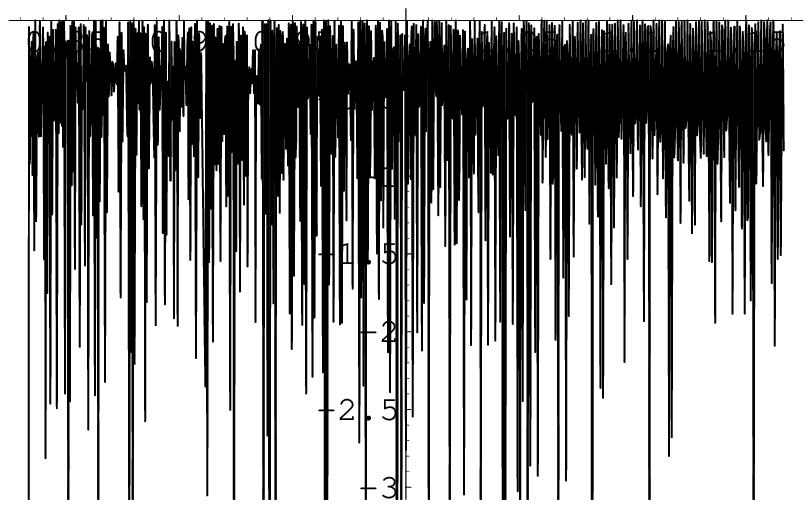}
\end{center}
\caption{$N=1000$}
\end{figure}

We also draw the graph of
\begin{equation*}
  \Re
  \left\{
     \log J_N\bigl(E;\exp(2\pi\sqrt{-1}/N)\bigr)
  \right\}
  -
  \Re
  \left\{
    \frac{N}{\sqrt{-1}}S(0)
    +\frac{3}{2}\log{N}
    +\frac{1}{2}\log\frac{T_E(0)}{2\pi^2}
  \right\}
\end{equation*}
for $N=100\times n$ with $1\le n\le 100$.
\begin{figure}[h]
\begin{center}
\includegraphics{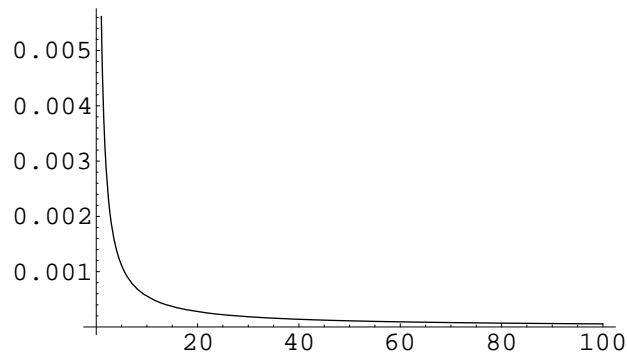}
\end{center}
\caption{The horizontal axis corresponds to $N/100$.}
\end{figure}
\bibliography{mrabbrev,hitoshi}
\bibliographystyle{hamsplain}
\end{document}